\documentclass[a4paper]{amsart}

\usepackage{amsmath,amsthm,amssymb,amscd}
\usepackage[arrow,matrix]{xy}
\usepackage[dvips]{graphics}
\usepackage{enumerate}

\theoremstyle{plain}
\numberwithin{equation}{section}
\newtheorem{thm}{Theorem}[section]
\newtheorem{prop}[thm]{Proposition}
\newtheorem{cor}[thm]{Corollary}
\newtheorem{lem}[thm]{Lemma}
\theoremstyle{definition}
\newtheorem{dfn}[thm]{Definition}
\newtheorem{exm}[thm]{Example}

\newtheorem{rem}[thm]{Remark}

\def\End{\mathop{\mathrm{End}}\nolimits}

\def\gr{\mathop{\mathrm{Gr}}\nolimits}

\def\CC{\mathbb{C}}
\def\QQ{\mathbb{Q}}
\def\RR{\mathbb{R}}
\def\ZZ{\mathbb{Z}}

\def\Aut{\mathop{\mathrm{Aut}}\nolimits}
\def\Ad{\mathop{\mathrm{Ad}}\nolimits}
\def\Ext{\mathop{\mathrm{Ext}}\nolimits}
\def\Ker{\mathop{\mathrm{Ker}}\nolimits}
\def\Gr{\mathop{\mathrm{Gr}}\nolimits}
\def\Spec{\mathop{\mathrm{Spec}}\nolimits}
\def\toric{\mathrm{toric}}
\def\torus{\mathrm{torus}}
\def\Gm{\mathbb{G}_{m}}
\def\KNU{J^{\mathrm{KNU}}}
\def\GGK{J^{\mathrm{GGK}}}
\def\GGKi{J^{\mathrm{GGK},0}}
\def\Flim{F_{\tilde{\phi}}}
\def\Flimp{F_{\tilde{\phi}}}
\def\Fnu{F_{\tilde{\nu}}}

\def\Fsigmap{F_{(\sigma' ,Z),\tilde{\phi}}}

\def\nusigmah{\hat{\nu}_{(\sigma',Z),\tilde{\phi}}}
\def\Dsigma{\Gamma(\sigma)^{\mathrm{gp}}\backslash D_{\sigma}}
\def\Dsigmap{\Gamma'(\sigma')^{\mathrm{gp}}\backslash D'_{\sigma'}}

\def\DSigmap{\Gamma' \backslash D'_{\Sigma'}}

\def\Esigma{E_{\sigma}}
\def\Esigmap{E'_{\sigma'}}
\def\Echeck{\check{E}_{\sigma}}
\def\Echeckp{\check{E}'_{\sigma'}}
\def\Jsigmap{J_{\sigma'}}
\def\Ksigmap{K_{\sigma'}}

\title[N\'eron models]{N\'eron models of Green-Griffiths-Kerr\\ and log N\'eron models}
\author[T.~Hayama]{Tatsuki HAYAMA}
\thanks{Supported by Grant-in-Aid for JSPS Fellows from Japan Society for the Promotion of Science.}
\address{Department of mathematics, National Taiwan University, Taipei 106, Taiwan}
\email{t-hayama@math.ntu.edu,tw}
\subjclass{14D07.} 
\keywords{N\'eron model; log mixed Hodge structure; admissible normal function; intermediate Jacobian}

\begin{document}
\maketitle
\begin{abstract}
For a variation of Hodge structure over a punctured disk, Green, Griffiths and Kerr introduced a N\'eron model which is a Hausdorff space that includes values of admissible normal functions.
On the other hand, Kato, Nakayama and Usui introduced a N\'eron model as a logarithmic manifold using log mixed Hodge theory.
This work constructs a homeomorphism between these two models.
\end{abstract}
\section{Introduction}
Let $J\to\Delta^*$ be a family of intermediate Jacobians arising from a variation of polarized Hodge structure (VHS) of weight $-1$ with a unipotent monodromy on a punctured disk.
By Carlson \cite{C}, the intermediate Jacobians are isomorphic to the extension groups of the Hodge structures,  in the category of mixed Hodge structures (MHS).
Then a section of $J\to\Delta^*$ is known as a variation of MHS (VMHS).
A VMHS satisfying the {\it admissibility condition} \cite{SZ} is called an admissible VMHS (AVMHS) and a section which gives an AVMHS is known as an admissible normal function (ANF) \cite{Sa1}.  

For the VHS, Green, Griffiths and Kerr \cite{GGK1} introduced the family $\GGK\to\Delta$ satisfying the following conditions:
\begin{itemize}
\item The family restricted to $\Delta^*$ is $J\to\Delta^*$;
\item The fiber over $0$ is a complex Lie group;
\item Any ANF is a section of $\GGK\to\Delta$;
\item $\GGK$ is a Hausdorff space.  
\end{itemize}
The total space $\GGK$ is called a {\it  N\'eron model}.
Here, $\GGK$ is simply a topological space.
The authors of \cite{GGK1} propose that ``One may `do geometry'" on N\'eron models.

In contrast, Kato, Nakayama and Usui constructed N\'eron models via a log mixed Hodge theory.
To explain their work, we describe $J\to\Delta^*$ by another formulation. 
Let $\Delta^*\to \Gamma\backslash D$ be the period map arising from the VHS.
The family of intermediate Jacobians can then be  obtained as the fiber product: 
$$\xymatrix{
J\ar@{->}[d]\ar@{->}[r]&\Gamma '\backslash D' \ar@{->}[d]^{\Gr^W_{-1}}\\
\Delta^*\ar@{->}[r]&\Gamma\backslash D
}
$$
where $D'$ and $\Gamma'$ are used for the MHS corresponding to the intermediate Jacobians.

Kato, Nakayama and Usui \cite{KNU1} extended the above diagram. 
First, Kato and Usui \cite{KU} stated that the period map can be extended to
\begin{align*}
\xymatrix{\Delta\ar@{}[d]|{\rotatebox{90}{$\subset$}}\ar@{->}[r]&\Gamma \backslash D_{\Sigma} \ar@{}[d]|{\rotatebox{90}{$\subset$}}\\
\Delta^*\ar@{->}[r]&\Gamma \backslash D}
\end{align*}
where $\Sigma$ is the fan of nilpotent cones arising from the monodromy of the VHS.
Here a boundary point of $\Gamma \backslash D_{\Sigma}$ is a {\it nilpotent orbit}, which approximates the period map given by Schmid \cite{Sc}.
The main theorem of \cite{KU} states that $\Gamma \backslash D_{\Sigma}$ is a logarithmic manifold and that it is a moduli space of log (pure) Hodge structures. 

Next, an ANF is written as
$$\Delta^*\to\Gamma '\backslash D'.$$
Kato, Nakayama and Usui \cite{KNU2} gives the fan $\Sigma'$, by which this map can be extended to 
\begin{align*}
\xymatrix{\Delta\ar@{}[d]|{\rotatebox{90}{$\subset$}}\ar@{->}[r]&\Gamma' \backslash D'_{\Sigma'} \ar@{}[d]|{\rotatebox{90}{$\subset$}}\\
\Delta^*\ar@{->}[r]&\Gamma' \backslash D'}.
\end{align*}
Similarly in the pure case \cite{KU}, a boundary point of $\Gamma' \backslash D'_{\Sigma'}$ is a {\it nilpotent orbit}, which approximates the ANF by the method proposed by Pearlstein \cite{P}.
The main theorem of \cite{KNU2} states that $\Gamma ' \backslash D'_{\Sigma'}$ is a logarithmic manifold and a moduli space of log mixed Hodge structures. 

Finally, they define the {\it log N\'eron model} $\KNU$ as the fiber product
$$\xymatrix{
\KNU\ar@{->}[d]\ar@{->}[r]&\Gamma '\backslash D'_{\Sigma'} \ar@{->}[d]^{\Gr^W_{-1}}\\
\Delta\ar@{->}[r]&\Gamma\backslash D_{\Sigma}
}
$$
in the category of logarithmic manifolds.
We remark that $\KNU$ is not only a topological space but also has a geometric structure as a {\it logarithmic manifold}.

However, \cite{KNU1} does not show the relationship between $\GGK$ and $\KNU$.
In fact, \S 8.2 of \cite{KNU1} states that the relationship is apparently unknown between $\KNU$ and the N\'eron model constructed by Green, Griffiths and Kerr.
Our main aim is to solve this problem.
\begin{thm}[Theorem \ref{main}]
$\GGK$ is homeomorphic to $\KNU$.
\end{thm}
We explain the key of the proof. 
By using the liftings in Equations (\ref{ANF-lift}) and (\ref{lift2}), we construct the bijective map between them (in Proposition \ref{bij} of this paper).
In section 5, we show that this map is a homeomorphism.
The diagram in Equation (\ref{com_dgm4}) and the admissibility condition (Equations (\ref{ANFdef}) or (\ref{fanp})) play important roles in the proof.

\subsection*{Acknowledgment}
The author is grateful to Professors Kazuya Kato, Chikara Nakayama, Gregory Pearlstein, Christian Schnell and Sampei Usui for their valuable advice and warm encouragement. 
\section{Preliminary}
In this section, we recall the definitions of the N\'eron models given in \cite{GGK1} and \cite{KNU2}.
Let $(\mathcal{H}_{\ZZ},\mathcal{F},\nabla )$ be a variation of polarized Hodge structure of weight $-1$ over a punctured disk $\Delta^*$, where $\mathcal{H}_{\ZZ}$ is a local system, $\mathcal{F}$ is a filtration of a locally free sheaf  $\mathcal{H}:=\mathcal{H}_{\ZZ}\otimes\mathcal{O}_{\Delta^*}$ and $\nabla$ is a Gauss-Manin connection.
We assume that the monodromy transformation $T$ is unipotent.
\subsection{Families of intermediate Jacobians}
Let $(H,F)$ be the total space of the vector bundle corresponding to the VHS $(\mathcal{H},\mathcal{F})$.
The intermediate Jacobian over $s\in\Delta^*$ is defined as
$$J_s:=F_s^0\backslash H_s/\mathcal{H}_{\ZZ;s}$$
where the subscript $s$ denotes the fiber (or stalk) over $s$.
By Carlson \cite{C}, we have the isomorphism
\begin{align}\label{ExtMHS}
\Ext_{\mathrm{MHS}}^1{(\ZZ(0),H_{s})}\cong J_{s}
\end{align}
where $\ZZ(0)$ is the Tate's Hodge structure.

We describe the family of intermediate Jacobians $J\to\Delta^*$ using the MHS in (\ref{ExtMHS}).
Fix a reference point $s_0\in\Delta^*$.
For the PHS $H_{s_0}=(H_{\ZZ},F_{s_0},\langle\; ,\;\rangle )$ over $s_0$, we take a MHS $H'$ which represents an extension class in $\Ext_{\mathrm{MHS}}^1{(\ZZ(0),H_{s_0})}$.ontinuation
Let $D$ (resp. $D'$) be the period domain for the type of $H_{s_0}$ (resp. $H'$), defined in \cite{G} (resp. \cite{U}).
The VHS gives the period map $\phi:\Delta^*\to\Gamma\backslash D$ where $\Gamma$ is the monodromy group.
Here we may write 
$$\Gamma=\{T^{n}\in\Aut{(H_{\ZZ})}\;|\; n\in\ZZ\}.$$
Then the family of intermediate Jacobians is obtained by the following Cartesian diagram:
$$\xymatrix{
J\ar@{->}[d]\ar@{->}[r]&\Gamma '\backslash D' \ar@{->}[d]^{\Gr^W_{-1}}\\
\Delta^*\ar@{->}[r]^{\phi}&\Gamma\backslash D
}
$$
where $\Gamma ':=\{T'\in \Aut{(H'_{\ZZ})}\;|\; T'|_{\Aut{(H_{\ZZ})}}\in \Gamma\}$.

We now review some properties of the period domains $D$ and $D'$.
Let $\check{D}$ (resp. $\check{D'}$) be the compact dual of $D$ (resp. $D'$), defined in \cite{G} (resp. \cite{U}).
From \cite[\S 4]{G} (resp. \cite[\S 2]{U}) we have the following properties for the pure case (resp. for some mixed case including the case of $D'$): 
\begin{prop}
Let $G_{A}:=\Aut{(H_{A},\langle\; ,\;\rangle)}$ (resp. $G'_{A}:=\Aut{(H'_{A},\langle\; ,\;\rangle_{\bullet })}$) for $A=\ZZ,\RR,\QQ,\CC$.
Then
\begin{enumerate}
\item $G_{\RR}$ (resp. $G'_{\RR}$) acts on $D$ (resp. $D'$) transitively;
\item $G_{\CC}$ (resp. $G'_{\CC}$) acts on $\check{D}$ (resp. $\check{D}'$) transitively;
\item Any subgroup of $G_{\ZZ}$ (resp. $G'_{\ZZ}$) acts on $D$ (resp. $D'$) properly discontinuously.
\end{enumerate}
\end{prop}

Since $H'$ is an extension of $H_{s_0}$ by $\ZZ(0)$, we have the exact sequence
$$0\to H_{\ZZ}\stackrel{i}{\to} H'_{\ZZ}\stackrel{j}{\to} \ZZ\to 0$$
of $\ZZ$-modules.
We fix $e\in H'_{\ZZ}$ such that $j(e)=1$.
Then we may write
\begin{align}\label{decomp}
H'_{\ZZ}\cong H_{\ZZ}\oplus \ZZ e.
\end{align}
We set
$$\mathfrak{h}:=\{X\in \End{(H'_{\CC})}\;|\; X|_{\End{(H_{\CC})}}=0,\;X(e)\in H_{\CC}\}.$$
\begin{prop}[\cite{U} Theorem 2.16]\label{fiberbdl}
$\Gr^{W}_{-1}:\check{D}'\to\check{D}$ is a fiber bundle with fiber  $\mathfrak{h}/(\mathfrak{h}\cap\mathfrak{b})$.
Here, $\mathfrak{b}$ is the Lie algebra of an isotropy subgroup of $G_{\CC}$.
\end{prop}

\subsection{Normal functions and the identity components}
We first define the following sheaves over $\Delta^*$:
\begin{align*}
&\mathcal{J}:=\mathcal{F}^0\backslash \mathcal{H}/\mathcal{H}_{\ZZ},\\
&\mathcal{J}_{\nabla}:=\left\{\begin{array}{l|l}\nu\in\mathcal{J} & \begin{array}{l}\nabla\tilde{\nu}\in \mathcal{F}^{-1}\otimes\Omega^1\\ \text{ for any local lifting }\tilde{\nu}\end{array}\end{array}\right\}.
\end{align*}
Since the monodromy is unipotent, we have the Deligne extension $(\mathcal{H}_e,\mathcal{F}_e)$.
Then we define the following sheaves over $\Delta$:
\begin{align*}
&\mathcal{J}_e:=\mathcal{F}^0_e\backslash \mathcal{H}_e/j_*\mathcal{H}_{\ZZ},\quad \mathcal{J}_{e,\nabla}:=\mathcal{J}_{e}\cap j_*\mathcal{J}_{\nabla}
\end{align*}
where $j:\Delta^*\hookrightarrow \Delta$.
A section of $\mathcal{J}_{e,\nabla}$ is called a {\it normal function} (NF).

Secondly we define a space that includes values of NF according to \cite[\S{\rm II}.A]{GGK1}.
Let $(H_e,F_e)$ be the total space of the vector bundles corresponding to $(\mathcal{H}_e,\mathcal{F}_e)$.
Since these vector bundles are trivial, we have a trivialization
\begin{align}\label{triv-can-ext}
F^n_e\cong\Delta\times F^n_{e;0}.
\end{align}
Since $(F_{e;0},W(N))$ is a MHS \cite{Sc}, we have the Deligne decomposition $H_{e;0}=\bigoplus_{p,q}I^{p,q}$.
This decomposition induces
\begin{align}\label{V}
F_{e;0}^0\backslash H_{e;0}\cong \bigoplus_{p<0}I^{p,q}=:V.
\end{align} 
By the trivialization (\ref{triv-can-ext}), we may write
$$F^0_e\backslash H_e\cong \Delta\times V.$$
We define the quotient space
$$J^Z:=F^0_e\backslash H_e/\sim$$
where the equivalence relation $\sim $ is given by equating 
two elements $(s,x),(s',x')\in \Delta\times V\cong F^0_e\backslash H_e$ if and only if $s=s'$ and $x-x'\in j_*\mathcal{H}_{\ZZ ;s}$.
We call it the Zucker space.

The Zucker space $J^Z$ includes values of NF.
However, $J^Z$ is {\it not} generally  a Hausdorff space (cf. \cite[{\rm II}.B.8]{GGK1}).
Hence, \cite{GGK1} defines the subspace of $J^Z$ so that it is a Hausdorff space including values of NF.

\begin{dfn}[\cite{GGK1} {\rm II}.A.9]\label{GGKtop}
Let
\begin{align}\label{W}
W:=\{(s,x)\in\Delta\times V\;|\; x\in\Ker{(N)} \text{ if }s=0\}.
\end{align}
The {\it identity component of the N\'eron model} is the subset $\GGKi :=W/\sim$ of the Zucker space $J^Z$.
Here the topology on $\GGKi$ is induced from the {\it strong topology} of $W$ in $\Delta\times V$ \cite[\S 3.1.1]{KU}.
\end{dfn}
The identity component has the following property:
\begin{prop}[\cite{GGK1} {\rm II}.A.9]
For a NF $\nu$, $\nu(0)\in \GGKi_0.$
\end{prop}
\begin{rem}
In \cite{GGK1}, the definition of the topology on $\GGKi$ seems to be unclear (a remark after \cite[Theorem {\rm II}.A.9]{GGK1} states ``This topology is modeled on the `strong topology' in \cite{KU}").
In this paper, we use the strong topology on $W\subset\Delta\times V$.
Saito \cite{Sa} shows the Hausdorff property in the case of the ordinary topology.
\end{rem}
\subsection{Admissible normal functions and N\'eron models}
In accord with \cite[\S {\rm II}.B]{GGK1}, we define the sheaf
\begin{align}\label{ANFdef}
\tilde{\mathcal{J}}_{e, \nabla}:=\left\{\begin{array}{l|l}\nu\in j_*\mathcal{J}_{\nabla}&\begin{array}{l}\tilde{\nu}\text{ has a logarithmic growth as a section of }\check{\mathcal{F}}_e^{0},\\(T-I)\tilde{\nu}\in (T-I)\mathcal{H}_{\QQ}\cap\mathcal{H}_{\ZZ}\text{ for any local lifting }\tilde{\nu}.\end{array}\end{array}\right\}
\end{align}
where we denote the analytic continuation around the origin $0$ of $\tilde{\nu}$ by $(T-I)\tilde{\nu}$.
A section of $\tilde{\mathcal{J}}_{e,\nabla}$ is called an {\it admissible normal function} (ANF).
By definition, we have the following exact sequence of sheaves:
\begin{align}\label{ANFseq}
0\to\mathcal{J}_{e,\nabla}\stackrel{i}{\to}\tilde{\mathcal{J}}_{e,\nabla}\stackrel{j}{\to} G_0\to 0.
\end{align}
Here $G_0$ is the skyscraper sheaf supported at $0$, whose stalk is
$$G:=\frac{(T-I)H_{\QQ}\cap H_{\ZZ}}{(T-I)H_{\ZZ}}.$$
We define the abelian group
$$\GGK_s:=\frac{\GGKi_s\times\tilde{\mathcal{J}}_{e,\nabla;s}}{\{(\nu(s),[i(\nu)]_s)\;|\;\nu\in\mathcal{J}_{e,\nabla}\}}$$
where $[i(\nu)]_s$ is the germ at $s\in\Delta$.
Since $\mathcal{J}_{e,\nabla;s}$ is a divisible abelian group (i.e., for every positive integer $n$ and every $\nu\in\mathcal{J}_{e,\nabla;s}$ there exists $\mu\in\mathcal{J}_{e,\nabla;s}$ such that $n\mu = \nu$) and $G$ is a finite group, the exact sequence of the stalks of (\ref{ANFseq}) is split \cite[{\rm II}.B.11]{GGK1}.
Then we have the isomorphism
$$\GGK_s\cong\begin{cases}\GGKi_s&\text{ if }s\neq 0,\cr\GGKi_s\times G&\text{ if }s=0.\end{cases}$$
\begin{dfn}[\cite{GGK1} {\rm II}.B.9]
The {\it N\'eron model} of Green-Griffiths-Kerr is the topological space
$$\GGK:=\bigsqcup_{s\in\Delta}\GGK_s.$$
Here the topology on $\GGK$ is defined by the open sets
\begin{align}\label{GGK_nbd}
S(\nu):=\left\{((s,x),[\nu]_{s})\in\GGK\;|\; (s,x)\in S\right\}
\end{align}
where $S$ is an open set of $\GGKi$ and $\nu$ is an ANF.
\end{dfn}
\begin{exm}[Classical case]
Let $\bar{f}:\bar{E}\to\Delta$ be a degenerating family of elliptic curves of Kodaira-type $I_n$.
For the restriction $f:E\to\Delta^*$, we have the local system $\mathcal{H}_{\ZZ}:=R^{1}f_*\ZZ$ and the filtration $\mathcal{F}^p=R^1f_*(\Omega_{E/\Delta^*}^{\geq p})$.
Here $(\mathcal{H}_{\ZZ},\mathcal{F})$ is a VHS over $\Delta^*$ with a unipotent monodromy.
In this case, 
$$\GGKi_0\cong\mathbb{G}_m,\quad G\cong \ZZ\backslash n\ZZ$$
twisting $(\mathcal{H}_{\ZZ},\mathcal{F})$ into the VHS of weight $-1$.
\end{exm}
\subsection{A nonclassical example}
We give an example where the N\'eron model is not an analytic space.
Our two sources, \cite[\S {\rm III}.A]{GGK2} and \cite[\S 9]{KNU1}, deal with special situations of this kind.

Let $Y$ be a singular $K3$ surface (i.e., $\rho (Y)=20$) and $\bar{f}:\bar{E}\to\Delta$ be a degenerating family of elliptic curves of Kodaira-type $I_n$.
By the Shioda-Inose correspondence \cite{SI}, for a transcendental basis $\{t_1,t_2\}$ of $H^2(Y)$, the intersection form is represented as 
$$(t_i\cdot t_j)_{i,j}=\begin{pmatrix} 2a &b\cr b &2c \end{pmatrix}$$
where $a,b,c\in \ZZ$, $a,c>0$ and $b^2-4ac<0$.
We assume that $a=m$, where $m$ is a square free positive integer, and that $b=0$ and $c=1$.
We take a symplectic basis $\{\alpha,\beta\}$ of $H^1(E_s)$ for $s\neq 0$ such that the monodromy action is 
$$\alpha\mapsto\alpha +n\beta,\quad \beta\mapsto\beta.$$
Setting 
$$e_1=t_1 \times \alpha,\quad e_2=t_2\times \alpha,\quad e_3=\frac{t_1}{2m}\times\beta,\quad e_4=\frac{t_2}{2}\times\beta$$
in $H^3(Y\times E_s,\QQ)$, then the intersection form is represented as 
$$(e_i\cdot e_j)_{i,j}=\begin{pmatrix} 0 &I\cr -I &0 \end{pmatrix}.$$

For the family $g:=f\circ \mathrm{pr}_2:Y\times E\to \Delta^*$, we set the local system $\mathcal{H}_{\ZZ}\subset R^{3}g_*\QQ$ such that $\mathcal{H}_{\ZZ,s}=\sum_i\ZZ e_i$ and the filtration $\mathcal{F}^p$ induced from $R^3 g_*(\Omega_{Y\times E/\Delta^*}^{\geq p})$.
Then $(\mathcal{H}_{\ZZ},\mathcal{F})$ is a VHS and a fiber $(\mathcal{H}_{\ZZ,s},F_s)$ is a PHS of weight $-1$ where $h^{1,-2}=h^{0,-1}=h^{-1,0}=h^{-2,1}=1$, twisting it into the VHS of weight $-1$.
The monodromy transformation is written by
\begin{equation*}
T=
\begin{pmatrix}\LARGE{I_2}&\LARGE{0}\cr \footnotesize{\begin{matrix}2mn&0\cr 0&2n\cr \end{matrix}}&\LARGE{I_2}\cr\end{pmatrix}.
\end{equation*}
By \cite[\S 12.3]{KU}, the limiting MHS is described by the following Hodge diamond:
$$\xymatrix{
\stackrel{(1,-1)}{\bullet}\ar@{->}[d]^N &\stackrel{(-1,1)}{\bullet}\ar@{->}[d]^N\\
\stackrel{(0,-2)}{\bullet} &\stackrel{(-2,0)}{\bullet}\\}$$
Then 
\begin{align*}
\GGKi_0&\cong I^{-2,0}/j_*\mathcal{H}_{\ZZ;0},\quad G\cong \ZZ/2mn\ZZ\times \ZZ/2n\ZZ.
\end{align*}
In this case, the dimension of $\GGKi_0$ is smaller than the dimension of a general fiber and $J^{\mathrm{Z}}$ is not a Hausdorff space (cf. \cite[\S 9]{KNU1}).

%


\subsection{Moduli spaces of log Hodge structures and log N\'eron models}
Let $\mathfrak{g}_{A}$ (resp. $\mathfrak{g}'_{A}$) be the Lie algebra of $G_{A}$ (resp. $G'_{A}$) for $A=\RR,\CC$.
Writing $\sigma=\RR_{\geq 0}N$ in $\mathfrak{g}_{\RR}$ with $N=\log{(T)}$, we set the fan $\Sigma:=\{\{0\},\sigma\}$ and the set 
\begin{align}\label{nilporbit}
D_{\Sigma}=\{(\sigma ,Z)\;|\;\sigma\in \Sigma,\;Z=\exp{(\sigma_{\CC})}F \text{ is a }\sigma\text{-nilpotent orbit}\}.
\end{align}
By \cite{KU}, the period map $\phi:\Delta^*\to\Gamma\backslash D$ extends to the {\it log} period map $\phi:\Delta\to \Gamma\backslash D_{\Sigma}$.

Following \cite{KNU1}, we define the fan
\begin{align}\label{fanp}
\Sigma':=\left\{\begin{array}{l|l}\RR_{\geq 0}N'& \begin{array}{l}N'\in\End{H'_{\QQ}}, N'|_{\End{H_{\QQ}}}=N,\\ N'(e)=N(a)\text{ for some }a\in H_{\QQ} \text{ such that }(T-I)a\in H_{\ZZ}\end{array}\end{array}\right\}.
\end{align}
\begin{prop}
Let $\sigma'\in\Sigma '$.
Then there exists a generator $N'\in \End{H'_{\QQ}}$ of $\sigma'$ such that $\exp{(N')}\in \Gamma'$, and $\Ad{(\gamma)}\sigma'\in\Sigma '$ for $\gamma\in \Gamma'$.
Therefore $\Gamma '$ is strongly compatible with $\Sigma '$.
\end{prop}
\begin{proof}
By definition, a generator of $\sigma'$ is written by
\begin{align*}
\begin{pmatrix}N&Na\cr 0&0\cr\end{pmatrix}
\end{align*}
with respect to the decomposition (\ref{decomp}) for some $a\in H_{\QQ}$.
Moreover, we may write
$$\Gamma'=\left\{\begin{array}{l|l}\begin{pmatrix}T^n&b\cr 0&1\cr\end{pmatrix}&\begin{array}{l}b\in H_{\ZZ},\;n\in\ZZ \end{array}\end{array}\right\}.$$
Since $(T-I)a\in H_{\ZZ}$, we have
$$\exp{(N')}=\begin{pmatrix}T&(T-I)a\cr 0&1\cr\end{pmatrix}\in\Gamma'.$$

For $\gamma=\begin{pmatrix}T^n&b\cr 0&1\cr\end{pmatrix}\in\Gamma '$, we have
\begin{align}\label{Adg}
\Ad{(\gamma)}N'=\begin{pmatrix}N&N(T^na-b)\cr 0&0\cr\end{pmatrix}.
\end{align}
Since $(T-I)(T^na-b)\in H_{\ZZ}$, it follows that $\Ad{(\gamma)}N'\in\Sigma'$.
\end{proof}
Similarly in (\ref{nilporbit}), $D'_{\Sigma '}$ is defined as the set of nilpotent orbits \cite[\S 2.1.3]{KNU2}.
Using the above proposition, we define the action 
$$\Gamma '\times D'_{\Sigma'}\to D'_{\Sigma'};\;(\gamma,(\sigma',Z))\mapsto (\Ad{(\gamma)}\sigma',\gamma Z)$$
and the orbit space $\DSigmap$.

The geometric structure on $\DSigmap$ is defined in \cite[\S 2.2.2]{KNU2}.
For a nilpotent cone $\sigma'\in\Sigma '$, we set the monoid 
$$\Gamma'(\sigma'):=\Gamma'\cap\exp{(\sigma')}$$
and the toric variety
$$\toric_{\sigma'}:=\Spec{(\CC[\Gamma'(\sigma')^{\vee}])}_{\mathrm{an}}\cong\CC.$$
Moreover, we define the analytic space 
$$\check{E}'_{\sigma'}:=\toric_{\sigma'}\times\check{D}'$$
and the subspace
\begin{align*}
&E'_{\sigma'}=\left\{\begin{array}{l|l}\left( s,F\right)\in\check{E}'_{\sigma'}&
\begin{array}{l}\exp{(l(s)N')}F\in D'\text{ if $s\neq 0$,}\\
\exp{(\sigma'_{\CC})}F\text{ is a nilpotent orbit if $s=0$}
\end{array}
\end{array}\right\}\end{align*}
where $l(s)$ is a branch of $(2\pi i)^{-1}\log{(s)}$.
The topology on $E'_{\sigma'}$ is the strong topology in $\check{E}'_{\sigma'}$.
We then have the map
$$\Esigmap\stackrel{p'_1}{\to}\Dsigmap\stackrel{ p'_2}{\to}\DSigmap;\quad (s,F)\mapsto \begin{cases}(0,\exp{(l(s)N')}F)&\text{ if }s\neq 0,\\ (\sigma',\exp{(\sigma'_{\CC})}F)&\text{ if }s=0.\end{cases}$$
The geometric structure on $\DSigmap$ is induced from $\Esigmap$ locally through this map.
Moreover Kato, Nakayama and Usui announced the following theorem:
\begin{thm}[\cite{KNU2} Main Theorem]
Similarly in the pure case (\cite[Main Theorem]{KU}), the following holds: 
\begin{enumerate}
\item $\Esigmap,\Dsigmap$ and $\DSigmap$ are logarithmic manifolds;
\item $\Esigmap\to\Dsigmap$ is a $\sigma'_{\CC}$-torsor;
\item $\Dsigmap\to\DSigmap$ is locally an isomorphism;
\item $\DSigmap$ is a moduli space of log mixed Hodge structures.  
\end{enumerate}
\end{thm}
\begin{dfn}[\cite{KNU1} \S 7]
The {\it log N\'eron model} is the fiber product
\begin{align}\label{fbrpdt1}
\xymatrix{
\KNU\ar@{->}[d]\ar@{->}[r]&\DSigmap \ar@{->}[d]^{\Gr^{W}_{-1}}\\
\Delta\ar@{->}[r]^{\phi}&\Gamma\backslash D_{\Sigma}
}
\end{align}
in the category $\mathcal{B}(\log)$ \cite[3.2.4]{KU}.
\end{dfn}

We describe the topology on $\KNU$.
We now have the following diagram:
\begin{align}\label{fbrpdt2}
\xymatrix{
K_{\sigma'}\ar@{->}[d]\ar@{->}[r]&E'_{\sigma'} \ar@{->}[d]^{ p'_1}\\
J_{\sigma'}\ar@{->}[d]\ar@{->}[r]&\Dsigmap \ar@{->}[d]^{ p'_2}\\
\KNU\ar@{->}[r]&\DSigmap
}
\end{align}
where $K_{\sigma'}$ and $J_{\sigma'}$ are the fiber products in $\mathcal{B}(\log)$.
Here the topology on $K_{\sigma'}$ is the strong topology in $\Delta\times\check{E}'_{\sigma'}$.
The topological structures of $J_{\sigma'}$ (resp. $\KNU$) are induced from $K_{\sigma'}$ through the morphism $K_{\sigma'}\to J_{\sigma'}$ (resp. $K_{\sigma'}\to \KNU$).
\section{The relationship between $E_{\sigma}\to \Dsigma$ and $E'_{\sigma'}\to \Gamma'(\sigma')^{\mathrm{gp}}\backslash D'_{\sigma'}$}
The results of this section can be easily verified using \cite{KNU2}; however the details will be useful in later sections.
In the following section, we regard $\Esigma$ (resp. $E'_{\sigma'}$) as a topological space whose topology is the strong topology in $\check{E}_{\sigma}$ (resp. $\check{E}'_{\sigma'}$).
\subsection{$\sigma_{\CC}$-action on $\Esigma$ and $\sigma'_{\CC}$-action on $\Esigmap$ }
For $\sigma=\RR_{\geq 0}N\in\Sigma$, we define the algebraic torus
$$\torus_{\sigma}:=\Spec{(\CC[\Gamma(\sigma)^{\vee \mathrm{gp}}])}_{\mathrm{an}}\cong\Gm$$
and the toric variety
$$\torus_{\sigma}:=\Spec{(\CC[\Gamma(\sigma)^{\vee}])}_{\mathrm{an}}\cong\CC.$$
We then have the surjective map
$$\sigma_{\CC}\to\torus_{\sigma};\quad wN\mapsto \exp{(2\pi \sqrt{-1}w)},$$
which induces the action
\begin{align*}
\sigma_{\CC}\times\toric_{\sigma}\to\toric_{\sigma};\quad (wN,s)\mapsto \exp{(2\pi \sqrt{-1}w)}s.
\end{align*}
For $\sigma'=\RR_{\geq 0}N'\in\Sigma '$, the $\sigma'_{\CC}$-action on $\toric_{\sigma '}$ is defined similarly.

By the correspondence $N\leftrightarrow N'$ (resp. $\exp{(N)}\leftrightarrow \exp{(N')}$), we have the isomorphism $\sigma_{\CC}\cong\sigma'_{\CC}$ (resp. $\toric_{\sigma}\cong\toric_{\sigma'}$) and the following commutative diagram:
\begin{align}\label{com_dgm1}
\xymatrix{
\sigma'_{\CC}\times\toric_{\sigma'}\ar@{->}[r]\ar@{->}[d]^{\rotatebox{90}{$ \sim $}}&\toric_{\sigma'}\ar@{->}[d]^{\rotatebox{90}{$ \sim $}}\cr
\sigma_{\CC}\times\toric_{\sigma}\ar@{->}[r]&\toric_{\sigma}.
}\end{align}

Moreover we define the $\sigma_{\CC}$-action 
$$\sigma_{\CC}\times\Esigma\to\Esigma;\quad (wN,(s,F))\mapsto (\exp{(2\pi \sqrt{-1}w)}s,\exp{(-wN)}F),$$
and the $\sigma'_{\CC}$-action on $E'_{\sigma '}$ is defined similarly.
Setting
$$\gr ^{W}_{-1}:\check{E}'_{\sigma'}\to \check{E}_{\sigma};\quad (s,F)\mapsto (s,\gr^W_{-1}{(F)}),$$
the diagram (\ref{com_dgm1}) induces the following commutative diagram:
\begin{align}\label{com_dgm2}
\xymatrix{
\sigma'_{\CC}\times E'_{\sigma'}\ar@{->}[d]\ar@{->}[r]&E'_{\sigma '}\ar@{->}[d]\ar@{}[r]|{\subset}&\Echeckp\ar@{->}[d]^{\Gr^W_{-1}}\\
\sigma_{\CC}\times E_{\sigma}\ar@{->}[r]&E_{\sigma}\ar@{}[r]|{\subset}&\Echeck.
}
\end{align}
\subsection{The torsor property of $\Esigmap\to\Dsigmap$}
\begin{lem}
The action of $\sigma'_{\CC}$ on $E'_{\sigma'}$ is proper and free. 
\end{lem}
\begin{proof}
Since the lower horizontal action in (\ref{com_dgm2}) is free \cite[(7.2.9)]{KU}, the upper horizontal action in (\ref{com_dgm2}) is also free.

The $\sigma'_{\CC}$-action is proper if and only if the following condition is satisfied:
\begin{itemize}
\item For $x',y'\in E'_{\sigma'}$, sequences $\{x'_{\lambda}\}$ in $E'_{\sigma'}$ and $\{h'_{\lambda}\}$ in $\sigma'_{\CC}$ such that $x'_{\lambda}\to x'$ and $h'_{\lambda}x'_{\lambda}\to y'$, there exists $h'\in \sigma'_{\CC}$ such that $h'_{\lambda}\to h'$.\end{itemize}
We will now show that the above condition holds.
Taking $x',y',\{x'_{\lambda}\},\{h'_{\lambda}\}$ as above, we let
$$x:=\gr^{W}_{-1}{(x')},\quad y=\gr^{W}_{-1}{(y')}, \quad h_{\lambda}:=h'_{\lambda}|_{\End{H_{\QQ}}}.$$
Since the $\sigma_{\CC}$-action is proper \cite[(7.2.2)]{KU}, there exists $h\in \sigma_{\CC}$ such that $h_{\lambda}\to h$.
By the isomorphism $\sigma\cong \sigma'$, there exists $h'\in\sigma'_{\CC}$ such that $h=h'|_{\End{H_{\QQ}}}$ and $h'_{\lambda}\to h'$.     
\end{proof}
\begin{lem}[\cite{KU} Lemma 7.3.3]\label{torsor_condition}
Let $H$ be a topological group and $X$ be a topological space, and assume that we have an action $H\times X\to X$ which is proper and free.
Furthermore assume that the following condition is satisfied:\begin{itemize}
\item For $x\in X$, there exists a topological space $S$, a morphism $\iota :S\to X$ and an open neighborhood $U$ of $1$ in $H$ such that $U\times S\to X;\; (h,s)\mapsto h\iota (s)$ induces an isomorphism onto an open set of $X$.
\end{itemize}
Then $X\to H\backslash X$ is an $H$-torsor.
\end{lem}
\begin{prop}[\cite{KNU2} Theorem A.(2)]\label{Esigma_torsor}
The action of $\sigma'_{\CC}$ on $E'_{\sigma'}$ satisfies the condition of Lemma \ref{torsor_condition}. 
Then $E'_{\sigma'}\to \Dsigmap$ is a $\sigma'_{\CC}$-torsor.
\end{prop}
\begin{proof}
Since $\sigma'(s)_{\CC}\hookrightarrow T_{\check{D}'}(F)$ for $(s,F)\in E'_{\sigma'}$ (in this case $\sigma'(s)=\sigma'$ if $s=0$, and $\sigma'(s)=0$ otherwise), the proof is the same as for the pure case (\cite[(7.3.5)]{KU}).
\end{proof}
Since $ p_1:\Esigma\to\Dsigma$ (resp. $ p_1':\Esigmap\to\Dsigmap$) is a $\sigma_{\CC}$-torsor ($\sigma'_{\CC}$-torsor), the diagram (\ref{com_dgm2}) induces the following property:
\begin{cor}\label{cart}
The commutative diagram
\begin{align}\label{com_dgm3}
\xymatrix{
\Esigmap\ar@{->}[rr]^{p'_1}\ar@{->}[d]_{\Gr_{-1}^W}&&\Dsigmap\ar@{->}[d]\cr
E_{\sigma}\ar@{->}^{p_1}[rr]&&\Dsigma
}\end{align}
is Cartesian.
\end{cor}
\subsection{Limiting Hodge filtrations and liftings of the period map}
Let $\tilde{\phi}$ be a local lifting of the period map $\phi$.
Then we have the holomorphic map
\begin{align}\label{hat}
\hat{\phi}:\Delta^*\to \check{D};\quad s\mapsto \exp{(-l(s)N)}\tilde{\phi}(s).
\end{align}
We call this an {\it untwisted period map}.
By \cite{Sc}, this map is extended over $\Delta$.
We denote $\hat{\phi}(0)$ by $F_{\tilde{\phi}}$.
Remark that $F_{\tilde{\phi}}$ depends upon the choice of local lifting $\tilde{\phi}$.
The untwisted map $\hat{\phi}$ gives a lifting
$$\Delta\to E_{\sigma};\quad s\mapsto (s,\hat{\phi}(s))$$
of $\phi$.
This gives the following diagram:
\begin{align}\label{com_dgm4}
\xymatrix{
\check{E}'_{\sigma'}\ar@{}[r]|{\supset}\ar@{->}[d]_{\Gr^{W}_{-1}}&\Esigmap\ar@{->}[rr]^{p_1'}\ar@{->}[d]&&\Dsigmap\ar@{->}[d]\ar@{->}[rr]^{p'_2}&&\DSigmap\ar@{->}[d]\cr
\check{E}_{\sigma}\ar@{}[r]|{\supset}&E_{\sigma}\ar@{->}[rr]^{p_1}\ar@{<-}[d]_{(id,\hat{\phi})}&&\Gamma(\sigma)^{\mathrm{gp}}\backslash D_{\sigma}\ar@{}[rr]|{=}&&\Gamma\backslash D_{\Sigma}\cr
&\Delta\ar@{->}[rru]_{\phi}
}\end{align}
for $\sigma'\in \Sigma'$ such that $\sigma'\neq \{0\}$. 

For $(s,F)\in \check{E}'_{\sigma'}$ such that $\Gr^{W}_{-1}(F)=F_{\hat{\phi}(s)}$, we have the exact sequence
$$0\to F_{\hat{\phi}(s)}^p\to F^p\to \CC\to 0$$
if $p\leq 0$, and $F_{\hat{\phi}(s)}^p\cong F^p$ otherwise.
Then 
\begin{align}\label{fil}
F^p=\begin{cases}\CC (z,1)+F^p_{\hat{\phi}(s)}&\text{if }p\leq  0,\\
F^p_{\hat{\phi}(s)} &\text{if }p> 0\end{cases}
\end{align}
where $(z,1)\in H'_{\CC}$ is represented with respect to the decomposition (\ref{decomp}).
By the admissibility condition (\ref{fanp}), a generator of $\sigma'\in\Sigma'$ can be written as
$$N'=\begin{pmatrix}N&Na\cr 0&0\cr\end{pmatrix}$$
for some $a\in H_{\QQ}$.
\begin{prop}\label{fil_Es}
For $(s,F)\in \check{E}'_{\sigma'}$ such that $\Gr^{W}_{-1}(F)=F_{\hat{\phi}(s)}$,
$$(s,F)\in \Esigmap\Longleftrightarrow \begin{cases}z\in H_{\CC}&\text{if }s\neq 0,\\ z+a\in F^0_{\tilde{\phi}}+\Ker{(N)}&\text{if }s= 0\end{cases} $$
where $z\in H_{\CC}$ is in (\ref{fil}).
\end{prop}
\begin{proof}
If $s\neq 0$, then
\begin{align*}\Gr^{W}_{-1}\left(\exp{(l(s)N')}F\right)&=\exp{(l(s)N)}\hat{\phi}(s)=\tilde{\phi}(s)\in D
\end{align*}
for any $z\in H_{\CC}$.
Then $(s,F)\in \Esigmap$  for any $z\in H_{\CC}$.

If $s=0$, then 
$$N(z+a)\in F^{-1}_{\tilde{\phi}}$$
by the transversality condition for nilpotent orbits.
Since $(F_{\tilde{\phi}},W(N))$ is a MHS and $N$ is a $(-1,-1)$-morphism, $N(z+a)\in F^{-1}_{\tilde{\phi}}$ if $z+a\in F^0_{\tilde{\phi}}+\Ker{(N)}$.
\end{proof}
\section{A bijection}
In this section, we define a bijective map between $\KNU$ and $\GGK$ as sets.
\subsection{A map from $\GGK$ to $\KNU$}
Let $\nu$ be an ANF.
The ANF $\nu$ defines an AVMHS
$$\nu:\Delta^*\to\Gamma'\backslash D'.$$
Taking a local lifting $\tilde{\nu}$ of $\nu $, we have a local lifting $\tilde{\phi}=\Gr^W_{-1}(\tilde{\nu})$ of $\phi$.
Let $N'$ be the logarithm of monodromy of $\tilde{\nu}$.
Similarly in (\ref{hat}), we define
$$\hat{\nu}:\Delta^*\to\check{D}';\quad s\mapsto\exp{(-l(s)N')}\tilde{\nu}(s).$$
By the admissibility condition (\ref{ANFdef}), $\hat{\nu}$ extend over $\Delta$ and $\sigma'=\RR_{\geq 0}N$ is in $\Sigma '$.
We denote $\hat{\nu}(0)$ by $F_{\tilde{\nu}}$.
By \cite{P}, $(\sigma', F_{\tilde{\nu}})$ is a nilpotent orbit.
We have the commutative diagram
\begin{align}\label{ANF-lift}
\xymatrix{
&\check{D}'\ar@{->}[d]^{\Gr^W_{-1}}\\
\Delta\ar@{->}[ru]^{\hat{\nu}}\ar@{->}[r]^{\hat{\phi}}&\check{D},
}
\end{align}
that is $\hat{\nu}$ is a lifting of $\hat{\phi}$.

We fix $F_{\tilde{\nu}}$ as a reference point of $\check{D}'$. 
By Proposition \ref{fiberbdl}, the vertical morphism of the above diagram is the fiber bundle with fiber $\mathfrak{h}/(\mathfrak{h}\cap \mathfrak{b})$.
Recall that
\begin{align*}
&\mathfrak{h}=\{X\in \End{(H'_{\CC})}\;|\; X|_{\End{(H_{\CC})}}=0,\;X(e)\in H_{\CC}\},\\
&\mathfrak{h}\cap\mathfrak{b}=\{X\in\mathfrak{h}\;|\;X(e)\in F^0_{\tilde{\phi}}\},\\
& V=\bigoplus_{p<0}I^{p,q}\quad \left(\text{i.e., } F_{\tilde{\phi}}^0\oplus V= H_{\CC}\right),
\end{align*}
and then
\begin{align}\label{fiber}
\mathfrak{h}/(\mathfrak{h}\cap\mathfrak{b})\cong V;\quad X_{v}\leftrightarrow  v
\end{align}
where $X_v\in\mathfrak{h}$ such that $X_{v}(e)=v$.
Taking a boundary point $((0,\dot{v}),[\nu]_0)\in\GGK$
where
\begin{align}\label{vmod}
\dot{v}=v\mod{H_{\ZZ}\cap\Ker{N}}
\end{align}
for some $v\in V\cap\Ker{(N)}$, we define
$$\alpha((0,\dot{v}),[\nu]_0):=(0,(\sigma', \exp{(\sigma'_{\CC})}\exp{(X_{-v})}\Fnu)).$$
By Proposition \ref{fil_Es}, $\alpha((0,\dot{v}),[\nu]_0)$ is in $\KNU$. 
\begin{lem}
$\alpha((0,\dot{v}),[\nu]_0)$ is well-defined.
\end{lem} 
\begin{proof}
We show that $\alpha((0,\dot{v}),[\nu]_0)$ does not depend on the choice of $v$ of (\ref{vmod}), a lifting $\tilde{\nu}$ and a representative $((0,\dot{v}),[\nu]_0)$.

First, we take $x\in H_{\ZZ}\cap \Ker{(N)}$.
By (\ref{Adg}), this gives $\Ad{(\gamma_x)}N'=N'$ for
\begin{align*}
\gamma_x =\begin{pmatrix}I&x\cr 0&1\cr\end{pmatrix}\in \Gamma'.
\end{align*}

Then
\begin{align*}
(\sigma',\exp{(\sigma'_{\CC })}\exp{(X_{-v+x})}F_{\tilde{\nu}})&=(\sigma',\exp{(\sigma'_{\CC })}\gamma_x\exp{(X_{-v})}F_{\tilde{\nu}})\\
&=\gamma_x(\sigma',\exp{(\sigma'_{\CC })}\exp{(X_{-v})}F_{\tilde{\nu}}).
\end{align*}

Next, we take another lifting $\gamma\tilde{\nu}$ for $\gamma\in\Gamma'$.
The monodromy cone that arises from $\gamma\tilde{\nu}$ is $\Ad{(\gamma)}\sigma'$ and $F_{\gamma\tilde{\nu}}=\gamma F_{\tilde{\nu}}$.
Since $v\in \Ker{N}$, we have
$$\exp{(X_{-v})}\gamma=\gamma\exp{(X_{-v})}.$$
Then
\begin{align*}
(\Ad{(\gamma)}\sigma'_{\CC},\exp{(\Ad{(\gamma)}\sigma'_{\CC})}\exp{(X_{-v})} F_{\gamma\tilde{\nu}})&=\gamma(\sigma',\exp{(\sigma'_{\CC })}\exp{(X_{-v})}F_{\tilde{\nu}}).
\end{align*}

Finally, we take $((0,\dot{v}_1),[\nu_1]_0)\sim((0,\dot{v}_2),[\nu_2]_0)$
and let
$$F^p_{\tilde{\nu}_i}=\begin{cases}\CC(z_i,1)+F^p_{\tilde{\phi}}&\text{if }p\leq 0,\\ F^p_{\tilde{\phi}}& \text{if }p>0,\end{cases} \quad\text{for }i=1,2,$$
where $\tilde{\nu}_i$ are local liftings.
Let $\mu=\nu_1-\nu_2$.
Then there exists a local lifting $\tilde{\mu}$ such that
$$F^p_{\tilde{\mu}}=\begin{cases}\CC(z_1-z_2,1)+F^p_{\tilde{\phi}}&\text{if }p\leq 0,\\ F^p_{\tilde{\phi}}& \text{if }p>0.\end{cases} $$
Since $((0,\dot{v}_1),[\nu_1]_0)\sim((0,\dot{v}_2),[\nu_2]_0)$, $\mu$ is a NF such that $\mu(0)=\dot{v}_1-\dot{v}_2\in \GGKi_0$.
Then there exists $v_1,v_2\in\Ker{(N)}\cap V$ such that 
\begin{align*}
\dot{v}_i=v_i\mod{H_{\ZZ}\cap\Ker{N}}.
\end{align*}
and $z_1-z_2=v_1-v_2$.

On the other hand, the logarithm of the monodromy of $\tilde{\nu}_i$ is described by
\begin{align*}
\begin{pmatrix}N&Na_i\cr 0&0\cr\end{pmatrix}
\end{align*}
for some $a_i\in H_{\QQ}$.
Then the logarithm of the monodromy of $\tilde{\mu}$ is
\begin{align*}
\begin{pmatrix}N&N(a_1-a_2)\cr 0&0\cr\end{pmatrix}.
\end{align*}
Since $\mu$ is a NF, $(T-I)(a_1-a_2)\in (T-I)H_{\ZZ}$ by the exact sequence (\ref{ANFseq}).
Then $a_1-a_2\in H_{\ZZ}$.
Setting
\begin{align*}
\gamma_{a_1-a_2} =\begin{pmatrix}I&a_1-a_2\cr 0&1\cr\end{pmatrix},
\end{align*}
we have
\begin{align*}
\Ad{(\gamma_{a_1-a_2})}\begin{pmatrix}N&Na_2\cr 0&0\cr \end{pmatrix} =\begin{pmatrix}N&Na_1\cr 0&0\cr\end{pmatrix}
\end{align*}
by (\ref{Adg}).
Since $\alpha((0,\dot{v_2}),[\nu_2]_0)$ does not depend on the choice of lifting, we may take $\gamma_{a_1-a_2}\tilde{\nu}_2$ as a lifting of $\nu_2$.
The monodromy cone that arises from $\gamma_{a_1-a_2}\tilde{\nu}_2$ is $\sigma'$.
Then
\begin{align*}
(\sigma',\exp{(\sigma'_{\CC})}\exp{(X_{-v_2})}F_{\tilde{\nu}_2})
&=(\sigma',\exp{(\sigma'_{\CC })}\exp{(X_{-v_2})}\exp{(X_{v_2-v_1})}F_{\tilde{\nu}_1})\\
&=(\sigma',\exp{(\sigma'_{\CC })}\exp{(X_{-v_1})}F_{\tilde{\nu}_1}).
\end{align*} 
\end{proof}
Therefore, $\alpha$ defines a map
$$\alpha :\GGK\to\KNU$$
where the restriction $\alpha|_{J}$ is canonical.
\subsection{A map from $\KNU$ to $\GGK$}
Let $\tilde{\phi} $ be a lifting of $\phi$.
By Corollary \ref{cart}, for $(0,(\sigma',Z))\in J_{\sigma'}$, we have $(0,F)\in E'_{\sigma '}$ such that 
\begin{align}\label{Fsigma}
\Gr^{W}_{-1}(0,F)=(0,\Flim),\quad  p'_1(0,F)=(\sigma',Z).
\end{align}
We denote this filtration by $F_{(\sigma',Z),\tilde{\phi}}$.
\begin{lem}\label{g_Flim}
For $\gamma\in\Gamma'$ such that $\gamma |_{\Aut{H_{\ZZ}}}=T^n$, $\gamma \exp{((m-n)N')}\Fsigmap=F_{\gamma(\sigma',Z),T^m\tilde{\phi}}$.
\end{lem}
\begin{proof}
By Proposition \ref{fil_Es}, there exists $x\in\Ker{(N)}$ such that
$$F_{(\sigma',Z),\tilde{\phi}}^p=\begin{cases}\CC (x-a,1)+\Flimp^p&\text{if }p\leq  0,\\
\Flimp^p &\text{if }p> 0.\end{cases}$$
Writing $\gamma=\begin{pmatrix}T^n&b\cr\footnotesize{0}&\footnotesize{1}\end{pmatrix}$ for some $b\in H_{\ZZ}$,
we have
$$\gamma \exp{((m-n)N')}\Fsigmap^p=\begin{cases}\CC (T^mx-T^na+b,1)+F^p_{T^m\tilde{\phi}}&\text{if }p\leq  0,\\
F^p_{T^m\tilde{\phi}} &\text{if }p> 0.\end{cases}$$
Since $x\in\Ker{(N)}$, 
\begin{align}\label{Ker}
T^mx-T^na+b=x-(T^na-b).
\end{align}
By (\ref{Adg}) and Proposition \ref{fil_Es}, $(0, \gamma\exp{((m-n)N')}\Fsigmap)\in E'_{\Ad{(\gamma)\sigma'}}$, which satisfies 
\begin{align*}
&\Gr^{W}_{-1}(0, \gamma\exp{((m-n)N')}\Fsigmap)=(0,F_{T^m\tilde{\phi}}),\\
& p'_1(0, \gamma\exp{((m-n)N')}\Fsigmap)=\gamma(\sigma',Z).
\end{align*}
\end{proof}

Let $\hat{\phi}:\Delta\to\check{D}$ be the untwisted period map.
Since $\check{D}'\to\check{D}$ is a fiber bundle, there exists a lifting of $\hat{\phi}$:
\begin{align}\label{lift2}
\xymatrix{
&\check{D}'\ar@{->}[d]^{\Gr^W_{-1}}\\
\Delta\ar@{->}[ru]^{\nusigmah}\ar@{->}[r]^{\hat{\phi}}&\check{D}
}
\end{align}
such that $\nusigmah(0)=\Fsigmap$, shrinking $\Delta$ if necessary.
We then have a holomorphic map
$$\Delta^*\to\Gamma'\backslash D'; \quad s\mapsto p'_2\circ p'_1(s,\nusigmah(s)),$$
which defines an AVMHS, i.e., an ANF.
Denoting this ANF by $\nu_{(\sigma',Z),\tilde{\phi}}$, 
we define 
$$\beta(0,(\sigma',Z)):=((0,0),[\nu_{(\sigma',Z),\tilde{\phi}}]_0)\in \GGK.$$
\begin{lem}
$\beta(0,(\sigma',Z))$ is well-defined.
\end{lem} 
\begin{proof}
We show that $\beta(0,(\sigma',Z))$ does not depend on the choice of $\hat{\nu}_{(\sigma',Z),\tilde{\phi}}$, $(\sigma',Z)$ and $\tilde{\phi}$.

If we take liftings $\hat{\nu}_{(\sigma',Z),\tilde{\phi}}$ and $\hat{\nu}'_{(\sigma',Z),\tilde{\phi}}$ such that
\begin{align*}
\hat{\nu}_{(\sigma',Z),\tilde{\phi}}(0)=\hat{\nu}'_{(\sigma',Z),\tilde{\phi}}(0)=F_{(\sigma' ,Z),\tilde{\phi}},
\end{align*}
then $\mu:=\nu_{(\sigma',Z),\tilde{\phi}}-\nu'_{(\sigma',Z),\tilde{\phi}}$ is a NF and  $\mu(0)=0\in\GGKi_0$.
Then 
$$((0,0),[\nu_{(\sigma',Z),\tilde{\phi}}]_0)\sim((0,0),[\nu'_{(\sigma',Z),\tilde{\phi}}]_0).$$

Moreover, by Lemma \ref{g_Flim}, 
$$\gamma\exp{((m-n)N')}\hat{\nu}_{(\sigma',Z),\tilde{\phi}}(0)=F_{\gamma(\sigma',Z),T^m\tilde{\phi}}.$$
If we take $\hat{\nu}_{\gamma(\sigma',Z),T^m\tilde{\phi}}=\gamma\exp{((m-n)N')}\hat{\nu}_{(\sigma',Z),\tilde{\phi}}$ as a lifting of $T^{m}\hat{\phi}$,
then $\nu_{(\sigma',Z),\tilde{\phi}}=\nu_{\gamma(\sigma',Z),T^m\tilde{\phi}}$.
\end{proof}
Then $\beta$ defines a map
$$\beta :\KNU\to \GGK$$
 where the restriction $\beta|_{J}$ is canonical.
\begin{prop}\label{bij}
$\alpha=\beta^{-1}$ and $\beta=\alpha^{-1}$, i.e., $\GGK$ is bijective to $\KNU$.
\end{prop}
\begin{proof}
For $((0,\dot{v}),[\nu]_0)\in \GGK$, we set $(0,(\sigma',Z)):=\alpha((0,\dot{v}),[\nu]_0)$.
By making suitable choice of $\tilde{\nu}$, $\tilde{\phi}$ and $v$, we have $F_{(\sigma' ,Z),\tilde{\phi}}=\exp{(X_{-v})}\Fnu$.
Therefore $\mu(0)=\dot{v}$ for $\mu=\nu-\nu_{(\sigma' ,Z),\tilde{\phi}}$, which induces
$$((0,\dot{v}),[\nu]_0)\sim ((0,0),[\nu_{(\sigma ',Z),\tilde{\phi}}]_0)=\beta (0,(\sigma',Z)).$$

On the other hand, for $(0,(\sigma',Z))\in \KNU$, we set $((0,0),[\nu]_0):=\beta(0,(\sigma',Z))$.
By making suitable choice of $\tilde{\nu}$, $(\sigma',Z)$ and $\tilde{\phi}$, we have $\Fnu=F_{(\sigma',Z),\tilde{\phi}}$.
Therefore,
$$(0,(\sigma',Z))=(0,(\sigma',\exp{(\sigma'_{\CC})}F_{\tilde{\nu}}))=\alpha ((0,0),[\nu]_0).$$
\end{proof}
\section{A homeomorphism}
In this section, we show the following main theorem
\begin{thm}\label{main}
$\GGK$ is homeomorphic to $\KNU$.
\end{thm}

To show continuity, we describe an open neighborhood in $\KNU$.
We recall that the topology on $\KNU$ is induced from $K_{\sigma}$ through the following diagram:
\begin{align*}
\xymatrix{
K_{\sigma'}\ar@{->}[d]\ar@{->}[r]&E'_{\sigma'} \ar@{->}[d]^{ p'_1}\\
J_{\sigma'}\ar@{->}[d]\ar@{->}[r]&\Dsigmap \ar@{->}[d]^{ p'_2}\\
\KNU\ar@{->}[d]\ar@{->}[r]&\DSigmap \ar@{->}[d]^{\Gr^{W}_{-1}}\\
\Delta\ar@{->}[r]^{\phi}&\Gamma\backslash D_{\Sigma}.
}
\end{align*}
We describe an open neighborhood in $\KNU$ using the following steps:
\begin{enumerate}
\item[Step 1.] Describe an open neighborhood in $\Esigmap$;
\item[Step 2.] Describe an open neighborhood in $\Dsigmap$; and
\item[Step 3.] Describe an open neighborhood in $\KNU$.
\end{enumerate}
Open neighborhoods in $\GGK$ are described in (\ref{GGK_nbd}).
Comparing these, we show that the bijection constructed in the last section is continuous.
\subsection{Proof of the main theorem}
\subsubsection*{Setting:}
We take a boundary point $(0,(\sigma',Z))\in \KNU$.
Setting a lifting $\tilde{\phi}$ of the period map $\phi$, we have the untwisted period map $\hat{\phi}:\Delta\to\check{D}$.
Since $\sigma_{\CC}\hookrightarrow T_{\check{D}}(\Flim)$, we may take a $\CC$-subspace $B$ of $\mathfrak{g}_{\CC}$ such that $B\oplus \sigma_{\CC}\cong T_{\check{D}}(\Flim)$.
An open neighborhood at $\Flim$ in $\check{D}$ is described by
$$\{\exp{(a_1)}\exp{(a_2)}\Flim\;|\; a_1 \in U_1,\; a_2 \in U_2\}\cong U_1\times U_2$$
where $U_1$ (resp. $U_2$) is a sufficiently small open neighborhood of $0$ in $\sigma_{\CC}$ (resp. $B$).
We assume that the image of $\hat{\phi}$ is included in this open neighborhood shrinking $\Delta$ if necessary.
We put $\hat{\phi}(s)=(\hat{\phi}_1(s),\hat{\phi}_2(s))$, where $\hat{\phi}_1:\Delta\to \sigma_\CC\cong\CC$ is a holomorphic function such that $\hat{\phi}_1(0)=0$.
By using the coordinate $t=\exp{(2\pi\sqrt{-1}\hat{\phi}_1(s))}\cdot s$ on $\Delta$ the untwisted period map is
\begin{align*}
\hat{\phi}(t)=\exp(-l(t)N)\tilde{\phi}(t)&=\exp{\left(-\hat{\phi}_1(s)N-l(s)N\right)}\tilde{\phi}(t)\\
&=\exp{(-\hat{\phi}_1(s)N)}\hat{\phi}(s).
\end{align*}
Then $\hat{\phi}_1(t)=0$ for $\hat{\phi}(t)=(\hat{\phi}_1(t),\hat{\phi}_2(t))\in U_1\times U_2$.
It is significant that $\Flim$, $\Fnu$ and $F_{(\sigma',Z),\tilde{\phi}}$ do not depend on this coordinate change (i.e., the bijection $\alpha$ is independent). 

\subsubsection*{Step 1 and Step 2:}
In the pure case, neighborhoods in $\Esigma$ and in $\Dsigma$ are described in \cite[(7.3.5)]{KU}.
We describe neighborhoods in $\Esigmap$ and in $\Dsigmap$ in a similar way.
Now we have the point $(0,F_{(\sigma',Z),\tilde{\phi}})\in E'_{\sigma'}$ as described in (\ref{Fsigma}).
Since $\gr^{W}_{-1}:\check{D}'\to\check{D}$ is a fiber bundle whose fiber is $V$, we have a local trivialization
\begin{align}\label{triv}
(\gr^{W}_{-1})^{-1}(U_1\times U_2)\cong U_1\times U_2\times V.
\end{align}
Since $F_{(\sigma',Z),\tilde{\phi}}\in(\Gr^W_{-1})^{-1}(\Flim)$, we can assume that $(0,0,0)$ corresponds to $F_{(\sigma',Z),\tilde{\phi}}$.
Using this local trivialization, an open neighborhood at $(0,\Fsigmap)$ in $\Echeckp$ can be described by
$$\{(a_0,(a_1,a_2,v))\;|\; a_0\in U_0,\; a_1 \in U_1,\; a_2 \in U_2,\; v\in U_3\}$$
where $U_0$ (resp. $U_3$) is a sufficiently small open neighborhood of $0$ in $\toric_{\sigma '}$ (resp. $V$).
Let
\begin{align*}
A'=\{(a_0,(0,a_2,v))\;|\; a_0\in U_0,\; a_2 \in U_2,\;v\in U_3\},\quad S'=A'\cap \Esigmap.
\end{align*}
Using the diagram (\ref{com_dgm2}), the $\sigma'_{\CC}$-action defines an open inclusion map 
$$U_1\times S'\hookrightarrow\Esigmap.$$
This inclusion map induces the open inclusion map
$$\sigma'_{\CC}\times S'\hookrightarrow  E'_{\sigma'},$$
shrinking $S'$ if necessary.
Then $ p'_1(S')$ is an open set of $\Dsigmap$ and $ p'_1(S')\cong S'$.
Moreover, $ p'_2\circ p'_1(S')$ is an open neighborhood of $(\sigma',Z)$ in $\DSigmap$. 

\subsubsection*{Step 3:}
Since $p'_1(S')$ (resp. $p'_2\circ p'_1(S')$) is an open neighborhood in $\Dsigmap$ (resp. $\DSigmap$), $p'_1\left((\Delta\times S')\cap \Ksigmap \right)$ (resp. $p'_2\circ p'_1\left((\Delta\times S')\cap \Ksigmap \right)$) is an open neighborhood of $\Jsigmap$ (resp. $\KNU$). 
Moreover, since $p'_1(S')\cong S'$, 
$$p'_1\left((\Delta\times S')\cap \Ksigmap \right)\cong (\Delta\times S')\cap \Ksigmap .$$

We describe $(\Delta\times S')\cap \Ksigmap$ explicitly.
By the diagram (\ref{com_dgm4}), we have the following commutative diagram:
\begin{align}
\xymatrix{
K_{\sigma'}\ar@{->}[dd]\ar@{->}[rr]& &E'_{\sigma'} \ar@{->}[dl]_{\Gr_{-1}^{W}}\ar@{->}[d]^{p'_2\circ p'_1}\\
 &E_{\sigma} \ar@{->}[dr]^{ p_1}& \DSigmap\ar@{->}[d]^{\Gr_{-1}^{W}}\\
\Delta\ar@{->}[rr]^{\phi}& &\Gamma\backslash D_{\Sigma}.
}
\end{align}
Then, for $(t,\xi)\in\Delta\times \Esigmap$, $(t,\xi)\in \Ksigmap$ if, and only if, 
$$\phi (t)=\Gr_{-1}^{W}\circ p'_2\circ p'_1(\xi)=p_1\circ \Gr_{-1}^{W}(\xi).$$ 
\begin{lem}\label{inverseimage}
$\left(( p_1)^{-1}(\phi(t))\right)\cap \Gr^W_{-1}(S')=(t,\hat{\phi}(t)).$
\end{lem}
\begin{proof}
Since $p_1((t,\hat{\phi}(t)))=\phi(t)$
and $p_1$ is a $\sigma_{\CC}$-torsor, the fiber is
$$( p_1)^{-1}(\phi(t))=\sigma_{\CC}\cdot(t,\hat{\phi}(t))=\{(\exp{(2\pi\sqrt{-1}x)}t,\exp{(-x N)}\hat{\phi}(t))\;|\;x\in \CC\}.$$
The intersection with $U_0\times U_1\times U_2$ is
\begin{align*}&(U_0\times U_1\times U_2)\cap ( p_1)^{-1}(\phi(t))\\
&=\{(\exp(2\pi\sqrt{-1}a_1)t,-a_1,\hat{\phi}_2(t))\;|\;\exp(2\pi\sqrt{-1}a_1)t\in U_0, -a_1\in U_1\}.
\end{align*}
On the other hand, for $(a_0,0,a_2,v)\in S'$
$$\Gr^W_{-1}((a_0,0,a_2,v))=(a_0,0,a_2).$$
Then $(a_0,0,a_2)\in( p_1)^{-1}(\phi(t))$ if, and only if, $a_0= t$ and $a_2=\hat{\phi}_2(t)$.
\end{proof}



\begin{lem}
\begin{align}\label{nbd_S}
(\Delta\times S')\cap K_{\sigma'}&=\left\{\begin{array}{l|l}\left(t,(t,0,\hat{\phi}_2(t),v)\right)&\begin{array}{l} t\in U_0\cap \Delta, \\v\in \Ker (N)\cap U_3 \text{ if } t=0,\\ v\in  U_3 \text{ if } t\neq 0\end{array}\end{array}\right\}.
\end{align}
\end{lem}
\begin{proof}
By Lemma \ref{inverseimage}, for $(a_0,0,a_2,v)\in S'$
$$\phi(t)=p_1\circ\Gr_{W}^{-1}((a_0,0,a_2,v))\Rightarrow  a_0=t\text{ and }a_2=\hat{\phi}_2(t).$$
By Proposition \ref{fil_Es}, if $t\neq 0$, then $(t,0,\hat{\phi}_2(t),v)\in  S'$ for $v\in U_3$.
If $t=0$, since $(0,0,v)\in U_1\times U_2\times U_3$ corresponds to $\exp{(X_{v})}F_{(\sigma',Z),\tilde{\phi}}$, then  $(0,0,0,v)\in S'$ for $v\in F^0_{\tilde{\phi}}+\Ker{(N)}$.
Since $V\bigoplus F^{0}_{\tilde{\phi}}=H_{\CC}$ by definition (\ref{V}), $v\in \Ker{(N)}$.
\end{proof}

\subsubsection*{Homeomorphism.}
Let $S:=W\cap (\Delta\times U_3)$ where $W$ is in (\ref{W}) and $S$ is endowed with the strong topology in $\Delta\times U_3$.
Then $S$ is homeomorphic to (\ref{nbd_S}).
For the local trivialization (\ref{triv}), we get 
$$\hat{\nu}:\Delta\to U_1\times U_2 \times U_3 \subset \check{D}';\quad t\mapsto (0,\hat{\phi}_2(t),0).$$
Then we have an ANF
$$\nu:\Delta \to \DSigmap; \quad t\mapsto p_2'\circ p_1'(t,\hat{\nu}(t)).$$
Following (\ref{GGK_nbd}) we set a neighborhood
$$\dot{S}(\nu)=\{((t,-\dot{v}),[\nu]_t)| (t,\dot{v})\in \dot{S}\}$$
at $\alpha^{-1}(0,(\sigma',Z))=((0,0),[\nu]_0)$ in $\GGK$ where $\dot{S}$ is the image of $S$ in the quotient space $W/\sim$.
Then $\alpha(\dot{S}(\nu))$ is the image of (\ref{nbd_S}) through $p'_2\circ p'_1$, which is a neighborhood of $(0,(\sigma',Z))$.
In fact 
\begin{align*}
\alpha ((t,-\dot{v}),[\nu]_0)=&\begin{cases}(0,\exp{(l(t)N')}\exp{(X_{v})}\hat{\nu}(t))&\text{if }t\neq 0\\(\sigma',\exp{(\sigma'_{\CC})}\exp{(X_{v})}\hat{\nu}(0))&\text{if }t=0\end{cases}\\
&=p'_2\circ p'_1\left(t,(t,0,\hat{\phi}_2(t),v)\right).
\end{align*}

\cite[\S 3.1]{KU} gives a fundamental system of neighborhoods at $(0,0)$ in $S$.
This fundamental system of neighborhoods defines a fundamental system of neighborhoods at $((0,0),[\nu]_0)$ in $\GGK$, which goes to a fundamental system of neighborhoods at $(0,(\sigma',Z))$ in $\KNU$ through $\alpha$.
Therefore, $\alpha$ is a homeomorphism.

\end{document}